\documentclass[11pt,letterpaper,reqno]{amsart}
\usepackage{amsmath,amsthm,amsfonts,amssymb,mathtools,algpseudocode}
\usepackage{comment,bookmark,bm,color,float}

\hypersetup{pdfstartview={FitH}}

\newtheorem{Theorem}{{\bf Theorem}}[section]

\newtheorem{Algorithm}[Theorem]{{\bf Algorithm}}

\newtheorem{Lemma}[Theorem]{{\bf Lemma}}
\newtheorem{Remark}[Theorem]{{\bf Remark}}

\numberwithin{equation}{section}

\newcommand{\calA}{\mathcal{A}}
\newcommand{\calB}{\mathcal{B}}
\newcommand{\calC}{\mathcal{C}}
\newcommand{\calX}{\mathcal{X}}
\newcommand{\calY}{\mathcal{Y}}
\newcommand{\calS}{\mathcal{S}}

\newcommand{\RR}{\mathbb{R}}
\newcommand{\diag}{\emph{diag}}
\newcommand{\off}{\text{off}}

\newcommand{\Proj}{\emph{Proj}}

\begin{document}

\title[A Jacobi-type method for the approximate orthogonal tensor diagonalization]{Convergence of a Jacobi-type method for the approximate orthogonal tensor diagonalization}
\author{Erna Begovi\'{c}~Kova\v{c}}\thanks{Erna Begovi\'{c} Kova\v{c}, Faculty of Chemical Engineering and Technology, University of Zagreb, Maruli\'{c}ev trg 19, 10000 Zagreb, Croatia, \texttt{ebegovic@fkit.hr}}

\thanks{This work has been supported in part by Croatian Science Foundation under the project UIP-2019-04-5200.}
\date{\today}

\subjclass[2010]{15A69, 65F25, 65F99}
\keywords{Jacobi-type methods, convergence, tensor diagonalization, tensor decompositions, SVD}

\maketitle

\begin{abstract}
For a general third-order tensor $\mathcal{A}\in\mathbb{R}^{n\times n\times n}$ the paper studies two closely related problems, an SVD-like tensor decomposition and an (approximate) tensor diagonalization.
We develop a Jacobi-type algorithm that works on $2\times2\times2$ subtensors and, in each iteration, maximizes the sum of squares of its diagonal entries.
We show how the rotation angles are calculated and prove convergence of the algorithm.
Different initializations of the algorithm are discussed, as well as the special cases of symmetric and antisymmetric tensors.
The algorithm can be generalized to work on higher-order tensors.
\end{abstract}

\section{Introduction}

Singular value decomposition is arguably the most powerful tool of numerical linear algebra.
It is not surprising that, when compared to the matrix SVD, the tensor generalization is significantly more complicated, see e.g.\@ \cite{Hackbusch-book12,VVM12,Grasedyck10,DeLHOSVD}.
We study the SVD-like tensor decomposition in the Tucker format,
\begin{equation}\label{decomp}
\calA=\calS\times_1U_1\times_2U_2\cdots\times_dU_d,
\end{equation}
where $\calA$ and $\calS$ are tensors of order $d$ and $U_1,U_2,\ldots,U_d$ are orthogonal matrices. Here, the tensor $\calS$ mimics the diagonal matrix of singular values from the matrix SVD. It is well known that, in the tensor case, one cannot expect to obtain a diagonal core tensor $\calS$. Hence, our goal will be to get a decomposition~\eqref{decomp} where $\calS$ is ``as diagonal as possible''. This SVD-like tensor decomposition problem is closely related to the tensor diagonalization problem. It has many applications in signal processing, blind source separation, and independent component analysis~\cite{CoJuBook,DeLDeMV2001,ComonSP}.

Problem~\eqref{decomp} for tensors of order $d=3$ has been studied by Moravitz Martin and Van Loan~\cite{MMVL08}. In their paper the authors use a Jacobi-type method to solve the maximization problem stated in~\eqref{maximization} below. Their numerical results suggest convergence, although a convergence proof is not provided.
If the tensor $\calS$ from~\eqref{decomp} is a diagonal tensor, then $\calA$ can be diagonalized using orthogonal transformations. Since a general tensor cannot be diagonalized, we aim to achieve an approximate diagonalization.
A similar problem for symmetric tensors has been studied in a series of papers by Comon, Li and Usevich~\cite{LUC20,LUC19,LUC18} where a Jacobi-type method is also a method of choice.

In this paper we develop a Jacobi-type algorithm with the same idea as in~\cite{MMVL08}, to maximize the sum of squares of the diagonal, but the algorithm itself is different from the one in~\cite{MMVL08}. Moreover, we prove the convergence of our algorithm. Our convergence results are alongside those for the symmetric case from~\cite{LUC20,LUC19,LUC18}.
We are concerned with general tensors, that is, we do not assume any tensor structure, except in Section~\ref{sec:structured}, where we discuss several special cases.

One can observe the problem~\eqref{decomp} either as a minimization problem where the goal is to minimize the off-diagonal norm of $\calS$,
$$\off^2(\calS)=\|\calS\|_F^2-\|\diag(\calS)\|_F^2\rightarrow\min,$$
or as a maximization problem,
\begin{equation}\label{maximization}
\|\diag(\calS)\|_F^2\rightarrow\max,
\end{equation}
where the square of the Frobenius norm of diagonal entries of $\calS$ is maximized. We are going to work with the formulation~\eqref{maximization}.
We mainly focus on tensors of order $d=3$ and develop a block coordinate descent Jacobi-type algorithm for finding the decomposition
$$\calA=\calS\times_1U\times_2V\times_3W,$$
such that
$$\|\diag(\calS)\|_F^2=\sum_{i=1}^n\calS_{iii}^2$$
is maximized. We prove that the algorithm converges to a stationary point of the objective function.
As it will be explained later in the paper, the algorithm can easily be generalized to tensors of order $d>3$.

Our algorithm for an approximate tensor diagonalization can also be used for a low-rank tensor approximation. We can approximate $\calA$ by a rank-$r$ tensor $\tilde{\calA}$ in the following way. Starting from the decomposition~\eqref{decomp} we form a diagonal $r\times r\times\cdots\times r$ order $d$ tensor $\mathcal{D}$ such that the diagonal elements of $\mathcal{D}$ are $r$ diagonal elements of $\calS$ with the highest absolute values. Moreover, for $i=1,\ldots,r$, we take $U_{i,r}$ as columns of $U_i$ corresponding to the selected diagonal elements. Then, the low-rank approximation is obtained as
$$\tilde{\calA}=\mathcal{D}\times_1U_{1,r}\times_2U_{2,r}\cdots\times_dU_{d,r}.$$

In Section~\ref{sec:agm} we describe the problem and construct the algorithm for solving the maximization problem~\eqref{maximization}. We  prove the previously mentioned convergence results in Section~\ref{sec:cvg}, while in Section~\ref{sec:numerical} we provide several numerical examples. Moreover, in Section~\ref{sec:structured} we study the special cases of symmetric and antisymmetric tensors.

\section{Orthogonal tensor decomposition}\label{sec:agm}

\subsection{Preliminaries and notation}

We use the tensor notation from~\cite{KB09}, which is commonly used in the papers dealing with numerical algorithms for tensors.
Notation from~\cite{Landsberg-book11} is also commonly used in multilinear algebra, but somewhat less frequently in its numerical aspects.

Tensors  of order three or higher are denoted by calligraphic letters, e.g.\@ $\calX$.
Tensor \emph{fibers} are vectors obtained from a tensor by fixing all indices but one. For a third-order tensor, its fibers are columns, rows, and tubes.
The mode-$m$ \emph{matricization} of a tensor $\calX\in\RR^{n_1\times n_2\times\cdots\times n_d}$ is an $n_m\times(n_1\cdots n_{m-1}n_{m+1}\cdots n_d)$ matrix $X_{(m)}$ obtained by arranging mode-$m$ fibers of $\calX$ into columns of $X_{(m)}$. In this paper we mainly work with $3$rd order tensors. Thus, we will have $m=1,2,3$.

The mode-$m$ \emph{product} of a tensor $\calX\in\RR^{n_1\times n_2\times\cdots\times n_d}$ with a matrix $A\in\RR^{p\times n_m}$ is a tensor
$\calY\in\RR^{n_1\times\cdots\times n_{m-1}\times p\times n_{m+1}\times\cdots\times n_d}$,
$$\calY=\calX\times_m A, \quad \text{such that} \quad Y_{(m)}=A X_{(m)}.$$
Two important properties of the mode-$m$ product are
\begin{align}
\calX\times_mA\times_nB & =\calX\times_nB\times_mA, \quad m\neq n, \label{product1} \\
\calX\times_nA\times_nB & =\calX\times_n(BA). \label{product2}
\end{align}
The \emph{norm} of $\calX$ is a generalization of the matrix Frobenius norm. It is given by
$$\|\calX\|_F=\sqrt{\sum_{i_1=1}^{n_1}\sum_{i_2=1}^{n_2}\cdots\sum_{i_d=1}^{n_d} x_{i_1i_2\ldots i_d}^2}.$$
To lighten the notation throughout the paper we are going to write this norm simply as $\|\calX\|$.
The \emph{inner product} of two tensors $\calX,\calY\in\RR^{n_1\times n_2\times\cdots\times n_d}$ is given by
$$\langle\calX,\calY\rangle=\sum_{i_1=1}^{n_1}\sum_{i_2=1}^{n_2}\cdots\sum_{i_d=1}^{n_d} x_{i_1i_2\ldots i_d}y_{i_1i_2\ldots i_d}.$$
It is straightforward to check that $\langle\calX,\calX\rangle=\|\calX\|^2$.

The \emph{Tucker decomposition} is a decomposition of a tensor $\calX$ into a core tensor $\calS$ multiplied by a matrix in each mode,
\begin{equation}\label{tucker}
\calX=\mathcal{S}\times_1M_1\times_2M_2\times_3\cdots\times_dM_d.
\end{equation}

Tensor $\calX\in\mathbb{R}^{n\times n\times n}$ is \emph{diagonal} when $\calX_{ijk}\neq0$ only if $i=j=k$,
that is, if $\off(\calX)=0$.

\subsection{Problem description}

Let $\calA\in\mathbb{R}^{n\times n\times n}$. We are looking for an orthogonal Tucker decomposition
\begin{equation}\label{svd-like}
\calA=\calS\times_1U\times_2V\times_3W,
\end{equation}
where $U,V,W\in\mathbb{R}^{n\times n}$ are orthogonal matrices and $\calS\in\mathbb{R}^{n\times n\times n}$ is a core tensor such that
\begin{equation}\label{Smaximization}
\|\diag(\calS)\|^2=\sum_{i=1}^n\calS_{iii}^2\rightarrow\max.
\end{equation}
From relation~\eqref{svd-like} tensor $\calS$ can be expressed as
$$\calS=\calA\times_1U^T\times_2V^T\times_3W^T.$$
Hence, in order to solve the problem defined by~\eqref{svd-like} and~\eqref{Smaximization} for a given tensor $\calA$, we need to find orthogonal matrices $U,V,W$ that maximize the objective function
\begin{equation}\label{maxfunction}
f(U,V,W)=\|\diag(\calA\times_1U^T\times_2V^T\times_3W^T)\|^2\rightarrow\max.
\end{equation}
We do this using a Jacobi-type method with a block coordinate descent approach.

For the sake of simplicity, our analysis is restricted to equal-sized modes. However, with a few technical adjustments, the same algorithm can be constructed for $\calA\in\mathbb{R}^{n_1\times n_2\times n_3}$. Then, in~\eqref{svd-like} we have
$U\in\mathbb{R}^{n_1\times n_1}$, $V\in\mathbb{R}^{n_2\times n_2}$, $W\in\mathbb{R}^{n_3\times n_3}$, and $\calS\in\mathbb{R}^{n_1\times n_2\times n_3}$.

\subsection{Jacobi-type algorithm}

We now describe the Jacobi-type algorithm for solving the maximization problem defined by~\eqref{maxfunction}. This is an iterative algorithm. Its $k$th iteration has the form
\begin{equation}\label{kth}
\calA^{(k+1)}=\calA^{(k)}\times_1R_{U,k}^T\times_2R_{V,k}^T\times_3R_{W,k}^T, \quad k\geq0, \quad \calA^{(0)}=\calA,
\end{equation}
where $R_{U,k},R_{V,k},R_{W,k}$ are plane rotations with the following structure,
\begin{equation}\label{rotation}
R(i,j,\phi)=\left[
    \begin{array}{ccccc}
      I &  &  &  &  \\
       & \cos\phi &  & -\sin\phi &  \\
       &  & I &  &  \\
       & \sin\phi &  & \cos\phi &  \\
       &  &  &  & I \\
    \end{array}
  \right]
  \begin{array}{l}
     \\
     i \\
     \\
     j \\
     \\
     \end{array}.
\end{equation}
Index pair $(i,j)$ in the rotation matrix~\eqref{rotation} is called a \emph{pivot position}. The set of all possible pivot positions is
$\{(i,j) \ : \ 1\leq i<j\leq n\}$.
In the $k$-th step, matrices $R_{U,k},R_{V,k},R_{W,k}$ have the same pivot position $(i_k,j_k)$, while the rotation angle $\phi_k$ is, in general, different for each matrix.

Our algorithm uses a block coordinate descent approach. This means that each iteration consists of three microiterations where we hold two variables constant and vary the third one. We have
\begin{align}
\calB^{(k)} & =\calA^{(k)}\times_1R_{U,k}^T, \label{tmultiplication1} \\
\calC^{(k)} & =\calB^{(k)}\times_2R_{V,k}^T, \label{tmultiplication2} \\
\calA^{(k+1)} & =\calC^{(k)}\times_3R_{W,k}^T. \label{tmultiplication3}
\end{align}
Here, by $\calB^{(k)}$ and $\calC^{(k)}$ we denote the intermediate steps.
Of course, if we combine all three microiterations together, using the properties of mode-$m$ product, namely~\eqref{product1} and~\eqref{product2}, we get back to the iteration step~\eqref{kth}.

Let us see how the rotation angles in matrices $R_{U,k},R_{V,k},R_{W,k}$ are computed. For a fixed iteration step $k$ we observe a $2\times2\times2$ subproblem. Assume that $(i_k,j_k)=(p,q)$, $1\leq p<q\leq n$. A subtensor of $\calA$ corresponding to an index pair $(p,q)$ is denoted by $\hat{\calA}$ and we can write it as
$$\hat{\calA}(:,:,1)=\left[
                         \begin{array}{cc}
                           a_{ppp} & a_{pqp} \\
                           a_{qpp} & a_{qqp} \\
                         \end{array}
                       \right], \quad
\hat{\calA}(:,:,2)=\left[
                         \begin{array}{cc}
                           a_{ppq} & a_{pqq} \\
                           a_{qpq} & a_{qqq} \\
                         \end{array}
                       \right].
$$
Then, the corresponding $2\times2\times2$ subproblem is to find $2\times2$ rotations $\hat{R}_U,\hat{R}_V,\hat{R}_W$ such that
$$\|\diag(\hat{S})\|^2=\sigma_{ppp}^2+\sigma_{qqq}^2\rightarrow\max,$$
where
$$\hat{\calS}=\hat{\calA}\times_1\hat{R}_U^T\times_2\hat{R}_V^T\times_3\hat{R}_W^T,$$
and
$$\hat{\calS}(:,:,1)=\left[
                         \begin{array}{cc}
                           \sigma_{ppp} & \sigma_{pqp} \\
                           \sigma_{qpp} & \sigma_{qqp} \\
                         \end{array}
                       \right], \quad
\hat{\calS}(:,:,2)=\left[
                         \begin{array}{cc}
                           \sigma_{ppq} & \sigma_{pqq} \\
                           \sigma_{qpq} & \sigma_{qqq} \\
                         \end{array}
                       \right].
$$

Taking only microiteration~\eqref{tmultiplication1} we calculate rotation angles for matrix $\hat{R}_U$. We have
$$\left[\begin{array}{cccc}
             b_{ppp} & b_{pqp} & b_{ppq} & b_{pqq} \\
             b_{qpp} & b_{qqp} & b_{qpq} & b_{qqq} \\
           \end{array}
         \right]=\left[
    \begin{array}{cc}
      \cos\phi & \sin\phi \\
      -\sin\phi & \cos\phi \\
    \end{array}
  \right]\left[
           \begin{array}{cccc}
             a_{ppp} & a_{pqp} & a_{ppq} & a_{pqq} \\
             a_{qpp} & a_{qqp} & a_{qpq} & a_{qqq} \\
           \end{array}
         \right].
$$
The rotation angle $\phi$ is chosen to maximize the function
\begin{equation}\label{g1}
g_1(\phi) = b_{ppp}^2+b_{qqq}^2 =(a_{ppp}\cos\phi+a_{qpp}\sin\phi)^2+(-a_{pqq}\sin\phi+a_{qqq}\cos\phi)^2. 
\end{equation}
Such $\phi$ must satisfy relation $g_1'(\phi)=0$. Taking the derivative of $g_1$ we get
\begin{align*}
g_1'(\phi) & =2(\cos\phi^2-\sin\phi^2)(a_{ppp}a_{qpp}-a_{pqq}a_{qqq}) +2\cos\phi\sin\phi(a_{pqq}^2+a_{qpp}^2-a_{ppp}^2-a_{qqq}^2) \\
& =2\cos(2\phi)(a_{ppp}a_{qpp}-a_{pqq}a_{qqq}) +\sin(2\phi)(a_{pqq}^2+a_{qpp}^2-a_{ppp}^2-a_{qqq}^2) =0.
\end{align*}
Dividing this relation by $\cos(2\phi)$ we obtain
\begin{equation}\label{tan2phi1}
\tan(2\phi)=\frac{2(a_{ppp}a_{qpp}-a_{pqq}a_{qqq})}{a_{ppp}^2+a_{qqq}^2-a_{pqq}^2-a_{qpp}^2}.
\end{equation}
Similarly, we find the rotation angles for matrices $\hat{R}_V$ and $\hat{R}_W$ as
\begin{equation}\label{tan2phi2}
\tan(2\phi)=\frac{2(b_{ppp}b_{pqp}-b_{qpq}b_{qqq})}{b_{ppp}^2+b_{qqq}^2-b_{qpq}^2-b_{pqp}^2}
\end{equation}
and
\begin{equation}\label{tan2phi3}
\tan(2\phi)=\frac{2(c_{ppp}c_{ppq}-c_{qqp}c_{qqq})}{c_{ppp}^2+c_{qqq}^2-c_{ppq}^2-c_{qqp}^2},
\end{equation}
respectively.

In the relations~\eqref{tan2phi1}--\eqref{tan2phi3} it is possible that both the numerator and the denominator are equal to zero. If that happens for one of those relations, we can skip the rotation in the corresponding direction and move on to the next one.
If this is the case for all pairs, the algorithm will be terminated and it should be restarted with preconditioning. This will be explained in Section~\ref{sec:structured} for the case of antisymmetric tensors.

Rotation angles in $\hat{R}_U,\hat{R}_V,\hat{R}_W$ do not need to be calculated explicitly. We only need the sine and the cosine of the corresponding angles. However, once we have formulas for computing $\tan(2\phi)$, there is still a problem of calculating efficiently $\sin\phi$ and $\cos\phi$. We will show how it is done for the rotation in the first mode. The procedure is the same in other modes.
We go back to the relation~\eqref{tan2phi1}. Denote
$$\lambda=2(a_{ppp}a_{qpp}-a_{pqq}a_{qqq})\text{sign}(a_{ppp}^2+a_{qqq}^2-a_{pqq}^2-a_{qpp}^2),$$
$$\mu=\vert a_{ppp}^2+a_{qqq}^2-a_{pqq}^2-a_{qpp}^2\vert.$$
Moreover, we denote $t=\tan\phi$. Using the double-angle formula for tangent,
$$\tan(2\phi)=\frac{2t}{1-t^2},$$
relation~\eqref{tan2phi1} reads
$$\frac{2t}{1-t^2}=\frac{\lambda}{\mu}.$$
This is a quadratic equation in $t$,
$\lambda t^2+2\mu t-\lambda=0,$
with solutions
$t_1=\frac{-\mu+\sqrt{\mu^2+\lambda^2}}{\lambda}, t_2=\frac{-\mu-\sqrt{\mu^2+\lambda^2}}{\lambda}.$
Note that the equation for $t_1$ is numerically unstable because catastrophic cancellation may occur. Therefore, we multiply both numerator and denominator by $\mu+\sqrt{\mu^2+\lambda^2}$. That way we attain a numerically stable expression
$t_1=\frac{\lambda}{\mu+\sqrt{\mu^2+\lambda^2}}.$
Finally,
$$\cos\phi_i=\frac{1}{\sqrt{1+t_i^2}}, \quad \sin\phi_i=\frac{t_i}{\sqrt{1+t_i^2}}=t_i\cos\phi_i, \quad i=1,2.$$
We calculate both solutions and use the one that gives the bigger value of the function~\eqref{g1}.

The order in which we choose pivot pairs is called \emph{pivot strategy}. In our algorithm the pivot strategy is assumed to be cyclic. We choose an ordering of pairs $(i,j)$, $1\leq i<j\leq n$, which makes one cycle. Then we repeat the same cycle of pivot pairs until the convergence criterium is satisfied. Common examples of cyclic pivot strategies are row-wise and column-wise strategies with corresponding ordering of pivot pairs defined by
\begin{equation}\label{Or}
\mathcal{O}_r=(1,2),(1,3),\ldots,(1,n),(2,3),\ldots,(2,n),\ldots,(n-1,n)
\end{equation}
and
\begin{equation}\label{Oc}
\mathcal{O}_c=(1,2),(1,3),(2,3),\ldots,(1,n),(2,n),\ldots,(n-1,n),
\end{equation}
respectively. The convergence results from Section~\ref{sec:cvg} hold for any cyclic strategy.
Nevertheless, to ensure convergence of the algorithm, pivot pairs should satisfy an additional condition. We only take a pivot pair $(i,j)$ such that (at least) one of the following inequalities is true,
\begin{equation}\label{pivotcond}
\vert\langle\nabla_Q f,Q\dot{R}(i,j,0)\rangle\vert \geq\eta\|\nabla_Q f\|_2, \quad \text{for } Q=U,V,W,
\end{equation}
where $0<\eta\leq\frac{2}{n}$, $\dot{R}(i,j,0)$ denotes $\frac{\partial}{\partial\phi}R(i,j,\phi)\Big{\vert}_{\phi=0}$, $f$ is defined in~\eqref{maxfunction}, and the projected gradient $\nabla_Q f$ will be defined in Subsection~\ref{sec:grad}.
If $(i,j)$ does not satisfy any of the conditions~\eqref{pivotcond}, then it will be skipped and we move to the next pair in the cycle. It will be shown in Lemma~\ref{tm:lemma2} that for each inequality~\eqref{pivotcond} it is always possible to find an appropriate pivot pair.

In the $k$th step of the algorithm, when we have $\calA^{(k)}$, $U_k,V_k,W_k$, we first compute the sine and the cosine of the rotation angle in the rotation matrix $R_{U,k}$. We compute the auxiliary tensor $\calB^{(k)}$,
$$\calB^{(k)} = \calA^{(k)}\times_1R_{U,k}^T = \calA\times_1(R_{U,k}^TU_k^T)\times_2V_k^T\times_3W_k^T,$$
and
$$U_{k+1}=U_kR_{U,k}.$$
Then we repeat this procedure in the other modes. This is summarized in Algorithm~\ref{agm:jacobi}.

\begin{Algorithm}\label{agm:jacobi}
\hrule\vspace{1ex}
\emph{Jacobi-type algorithm for the approximate tensor diagonalization}
\vspace{0.5ex}\hrule
\begin{algorithmic}
\State \textbf{Input:} $\calA\in\mathbb{R}^{n\times n\times n}$.
\State \textbf{Output:} orthogonal matrices $U,V,W$
\State $k=0$
\State $\calA^{(0)}=\calA$
\State $U_0=V_0=W_0=I$
\Repeat
\State Choose pivot pair $(i,j)$.
\If {$(i,j)$ satisfies~\eqref{pivotcond} for $Q=U$}
\State Find $\cos\phi_k$ and $\sin\phi_k$ for $R_{U,k}$ using~\eqref{tan2phi1}.
\State $\calB=\calA^{(k)}\times_1R_{U,k}$
\State $U_{k+1}=U_kR_{U,k}$
\EndIf
\If {$(i,j)$ satisfies~\eqref{pivotcond} for $Q=V$}
\State Find $\cos\phi_k$ and $\sin\phi_k$ for $R_{V,k}$ using~\eqref{tan2phi2}.
\State $\calC=\calB\times_2R_{V,k}$
\State $V_{k+1}=V_kR_{V,k}$
\EndIf
\If {$(i,j)$ satisfies~\eqref{pivotcond} for $Q=W$}
\State Find $\cos\phi_k$ and $\sin\phi_k$ for $R_{W,k}$ using~\eqref{tan2phi3}.
\State $\calA^{(k+1)}=\calC\times_3R_{W,k}$
\State $W_{k+1}=W_kR_{W,k}$
\EndIf
\Until{convergence}
\end{algorithmic}
\hrule
\end{Algorithm}

We have several remarks regarding the Algorithm~\ref{agm:jacobi}.
\begin{itemize}
\item Algorithm~\ref{agm:jacobi} employs the identity initialization $U_0=V_0=W_0=I$. This is not necessarily done this way and it will be further discussed within the numerical examples in Section~\ref{sec:numerical}, as well as in relation with the antisymmetric tensors in Section~\ref{sec:structured}.
\item It is not needed to explicitly form rotation matrices and tensor matricizations in order to perform mode-$n$ multiplications in the algorithm.
\item Conditions on pivot pairs~\eqref{pivotcond} can be simplified to lower the computational effort. This will be shown after Lemma~\ref{tm:lemma2}. Moreover, the coefficient $0<\eta\leq\frac{2}{n}$ can vary, which will be examined in Section~\ref{sec:numerical}.
\end{itemize}

This algorithm can be generalized for the order-$d$ tensors where $d>3$. In that case we need to obtain orthogonal matrices $U_1,U_2,\ldots U_d$ such that maximization condition~\eqref{maximization} holds.
One iteration of the algorithm consists of $d$ microiterations that are analogues of those in~\eqref{tmultiplication1}, \eqref{tmultiplication2}, and~\eqref{tmultiplication3}.

\section{Convergence results}\label{sec:cvg}

\subsection{Gradient of the objective function}\label{sec:grad}

Before we move on to the convergence of the Algorithm~\ref{agm:jacobi}, let us say something about the gradient of the objective function $f:O_n\times O_n\times O_n\rightarrow\mathbb{R}$,
\begin{equation}\label{of}
f(U,V,W)=\|\diag(\calA\times_1U^T\times_2V^T\times_3W^T)\|^2,
\end{equation}
where $O_n$ stands for the group of orthogonal matrices of order $n$.
To calculate $\nabla f$ we need an auxiliary function
$\tilde{f}:\mathbb{R}^{n\times n}\times\mathbb{R}^{n\times n}\times\mathbb{R}^{n\times n}\rightarrow\mathbb{R}$ defined by the same formula~\eqref{of} as $f$. In other words, function $\tilde{f}$ is such that $f$ is the restriction of $\tilde{f}$ to the set of triples of orthogonal matrices. Then $\nabla f$ is the projection of $\nabla\tilde{f}$ onto the tangent space at $(U,V,W)$ to the manifold $O_n\times O_n\times O_n$. We have
\begin{align}
\nabla f(U,V,W) & =\left[
        \begin{array}{ccc}
          \nabla_U f(U,V,W) & \nabla_V f(U,V,W) & \nabla_W f(U,V,W) \\
        \end{array}
      \right] \nonumber \\
& =\Proj\left[
        \begin{array}{ccc}
          \nabla_U \tilde{f}(U,V,W) & \nabla_V \tilde{f}(U,V,W) & \nabla_W \tilde{f}(U,V,W) \\
        \end{array}
      \right] \nonumber \\
& =\left[\begin{array}{ccc}
U\Lambda(U) & V\Lambda(V) & W\Lambda(W) \\
\end{array}\right], \label{grad}
\end{align}
where
\begin{equation}\label{lambda}
\Lambda(Q):=\frac{Q^T\nabla_Q \tilde{f}-(\nabla_Q \tilde{f})^TQ}{2}.
\end{equation}

To calculate $\nabla f(U,V,W)$ we write $\tilde{f}$ as
$$\tilde{f}(U,V,W)=\|\diag(\calA\times_1U^T\times_2V^T\times_3W^T)\|^2=\sum_{l=1}^n\left(\sum_{i,j,k=1}^n a_{ijk}u_{il}v_{jl}w_{kl}\right)^2.$$
Element-wise, we get
\begin{align*}
\frac{\partial\tilde{f}}{\partial u_{ml}} & =2\left(\sum_{i,j,k=1}^n a_{ijk}u_{il}v_{jl}w_{kl}\right)
\left(\sum_{j,k=1}^n a_{mjk}v_{jl}w_{kl}\right) \\
& =2\left(\calA\times_1U^T\times_2V^T\times_3W^T\right)_{lll}\left(\calA\times_2V^T\times_3W^T\right)_{mll}, \\
\frac{\partial\tilde{f}}{\partial v_{ml}}
& =2\left(\calA\times_1U^T\times_2V^T\times_3W^T\right)_{lll}\left(\calA\times_1U^T\times_3W^T\right)_{lml}, \\
\frac{\partial\tilde{f}}{\partial w_{ml}}
& =2\left(\calA\times_1U^T\times_2V^T\times_3W^T\right)_{lll}\left(\calA\times_1U^T\times_2V^T\right)_{llm}.
\end{align*}
Then, we can use the above relations together with~\eqref{lambda} to get an explicit expression for $\Lambda(U)$,
\begin{align*}
(\Lambda(U))_{lp} & = \frac{1}{2}\left(\sum_{m=1}^n u_{ml}\frac{\partial\tilde{f}}{\partial u_{mp}}-\sum_{m=1}^n\frac{\partial\tilde{f}}{\partial u_{ml}}u_{mp}\right) \\
& = \left(\calA\times_1U^T\times_2V^T\times_3W^T\right)_{ppp}\left(\calA\times_1U^T\times_2V^T\times_3W^T\right)_{lpp} \\
& \qquad -\left(\calA\times_1U^T\times_2V^T\times_3W^T\right)_{lll}\left(\calA\times_1U^T\times_2V^T\times_3W^T\right)_{pll}.
\end{align*}
Similarly we get the expressions for $\Lambda(V)$ and $\Lambda(W)$.

The gradient of the objective function will be needed in order to prove the following convergence result and also to check the pivot conditions~\eqref{pivotcond}.

\subsection{Convergence theorem}

The convergence of Algorithm~\ref{agm:jacobi} is given in Theorem~\ref{tm:cvg}.

\begin{Theorem}\label{tm:cvg}
Every accumulation point $(U,V,W)$ obtained by Algorithm~\ref{agm:jacobi} is a stationary point of the function $f$ defined by~\eqref{of}.
\end{Theorem}

The proof follows the idea from~\cite{Ishteva13}, which was also used in~\cite{LUC18}, as well as in~\cite{Begovic21normal,BeKre17}. The major obstacle is that here we have a function of three variables, while earlier this procedure was used with single-variable functions. We prove Lemma~\ref{tm:lemma2}, which is an adaptation of Lemma 3.1 from~\cite{LUC18}. Then, using Lemma~\ref{tm:lemma2} we prove Lemma~\ref{tm:lemma3}, which is an essential step in the proof of Theorem~\ref{tm:cvg}.

\begin{Lemma}\label{tm:lemma2}
For any differentiable function $f:O_n\times O_n\times O_n\rightarrow\mathbb{R}$, $U,V,W\in O_n$, and $0<\eta\leq\frac{2}{n}$ it is always possible to find index pairs $(i_U,j_U)$, $(i_V,j_V)$, $(i_W,j_W)$ satisfying pivot condition~\eqref{pivotcond}.
\end{Lemma}

\begin{proof}
Observe that
$$\dot{R}(i,j,0)=e_je_i^T-e_ie_j^T.$$

From the definition of the operator $\Lambda$ we see that matrix $\Lambda(U)$ is skew-symmetric. Then, from the fact that the Euclidean norm is invariant under unitary transformations and from relation~\eqref{grad} we have
\begin{align}
\vert\langle\nabla_U f(U,V,W),U\dot{R}(i,j,0)\rangle\vert & =\vert\langle U\Lambda(U),U\dot{R}(i,j,0)\rangle\vert \nonumber \\
& =\vert\langle \Lambda(U),\dot{R}(i,j,0)\rangle\vert=2\vert\Lambda(U)_{ij}\vert. \label{lm:gradpom}
\end{align}
We can always find an index pair $(i_U,j_U)$ such that
$$\vert\Lambda(U)_{i_Uj_U}\vert\geq\frac{1}{n}\|\Lambda(U)\|_2.$$
Inserting this into equation~\eqref{lm:gradpom} with $(i,j)=(i_U,j_U)$, we get
$$\vert\langle\nabla_U f(U,V,W),U\dot{R}(i_U,j_U,0)\rangle\vert \geq\frac{2}{n}\|\Lambda(U)\|_2\geq\eta\|\Lambda(U)\|_2=\eta\|\nabla_U f(U,V,W)\|_2,$$
which proves assertion $(i)$.
Since the matrices $\Lambda(V)$ and $\Lambda(W)$ are also skew-symmetric, in the same way we obtain
\begin{align*}
\vert\langle\nabla_V f(U,V,W),V\dot{R}(i_V,j_V,0)\rangle\vert & \geq\frac{2}{n}\|\Lambda(V)\|_2\geq\eta\|\Lambda(V)\|_2=\eta\|\nabla_V f(U,V,W)\|_2, \\
\vert\langle\nabla_W f(U,V,W),W\dot{R}(i_W,j_W,0)\rangle\vert & \geq\frac{2}{n}\|\Lambda(W)\|_2\geq\eta\|\Lambda(W)\|_2=\eta\|\nabla_W f(U,V,W)\|_2.
\end{align*}
This proves assertions $(ii)$ and $(iii)$, respectively.
\end{proof}

\begin{Remark}
Conditions~\eqref{pivotcond} are equivalent to
$$2\vert\Lambda(Q)_{ij}\vert \geq\eta\|\Lambda(Q)\|_2, \quad \text{for } Q=U,V,W,$$
where $\Lambda(\cdot)$ is as in relation~\eqref{lambda}.
\end{Remark}

\begin{Lemma}\label{tm:lemma3}
Let $U_k$, $V_k$, $W_k$, $k\geq0$ be the sequences generated by Algorithm~\ref{agm:jacobi}.
Let $\overline{U},\overline{V},\overline{W}$ be a triple of orthogonal matrices satisfying $\nabla f(\overline{U},\overline{V},\overline{W})\neq0$.
Then there exist $\epsilon>0$ and $\delta>0$ such that
$$\|U_k-\overline{U}\|_2<\epsilon, \quad \|V_k-\overline{V}\|_2<\epsilon, \quad \|W_k-\overline{W}\|_2<\epsilon$$
implies
\begin{equation}\label{lm:fdelta}
f(U_{k+1},V_{k+1},W_{k+1})-f(U_k,V_k,W_k)\geq\delta.
\end{equation}
\end{Lemma}

\begin{proof}
Let us fix the iteration step $k$.
To shorten the notation set $\phi_U=\phi_{U,k}$, $\phi_V=\phi_{V,k}$, $\phi_W=\phi_{W,k}$, and
$R_{U,k}=R(i_k,j_k,\phi_U)$, $R_{V,k}=R(i_k,j_k,\phi_V)$, $R_{W,k}=R(i_k,j_k,\phi_W)$.
We define three functions $h_k^{(1)},h_k^{(2)},h_k^{(3)}:\mathbb{R}\rightarrow\mathbb{R}$,
\begin{align*}
h_k^{(1)}(\phi_1) & =f(U_kR(i_k,j_k,\phi_1),V_k,W_k), \\
h_k^{(2)}(\phi_2) & =f(U_kR_{U,k},V_kR(i_k,j_k,\phi_2),W_k), \\
h_k^{(3)}(\phi_3) & =f(U_kR_{U,k},V_kR_{V,k},W_kR(i_k,j_k,\phi_3)).
\end{align*}
Further on, we define yet another function $h_k:\mathbb{R}\times\mathbb{R}\times\mathbb{R}\rightarrow\mathbb{R}$,
$$h_k(\phi_1,\phi_2,\phi_3)=f(U_kR(i_k,j_k,\phi_1),V_kR(i_k,j_k,\phi_2),W_kR(i_k,j_k,\phi_3)).$$
Since $R(i_k,j_k,0)=I$,
$$h_k(0,0,0)=f(U_k,V_k,W_k).$$
From the construction of Algorithm~\ref{agm:jacobi} we know that
\begin{align*}
\max_{\phi_1}h_k^{(1)}(\phi_1) & =h_k^{(1)}(\phi_U)=f(U_kR_{U,k},V_k,W_k), \\
\max_{\phi_2}h_k^{(2)}(\phi_2) & =h_k^{(2)}(\phi_V)=f(U_kR_{U,k},V_kR_{V,k},W_k), \\
\max_{\phi_3}h_k^{(3)}(\phi_3) & =h_k^{(3)}(\phi_W)=f(U_kR_{U,k},V_kR_{V,k},W_kR_{W,k}),
\end{align*}
and the $k$th step of the algorithm is represented by
$$h_k(\phi_U,\phi_V,\phi_W)=f(U_kR_{U,k},V_kR_{V,k},W_kR_{W,k})=f(U_{k+1},V_{k+1},W_{k+1}).$$
Moreover, it is easy to see from the algorithm that
\begin{equation}\label{lm:fineq}
f(U_{k+1},V_{k+1},W_{k+1}) \geq f(U_{k+1},V_{k+1},W_k) \geq f(U_{k+1},V_k,W_k) \geq f(U_k,V_k,W_k).
\end{equation}
In order to attain inequality~\eqref{lm:fdelta} we need at least one sharp inequality in~\eqref{lm:fineq}.

If $\nabla f(\overline{U},\overline{V},\overline{W})\neq0$, then at least one partial gradient of $f$ is not zero, that is
$$\nabla_U f(\overline{U},\overline{V},\overline{W})\neq0, \ \nabla_V f(\overline{U},\overline{V},\overline{W})\neq0,
\text{ or } \nabla_W f(\overline{U},\overline{V},\overline{W})\neq0.$$
Let us assume that $\nabla_U f(\overline{U},\overline{V},\overline{W})\neq0$. Then there exists $\epsilon>0$ such that
\begin{equation}\label{lm:mu1}
\mu_1:=\min\{\|\nabla_U f(U,V,W)\|_2 \ : \ \|U-\overline{U}\|_2<\epsilon\}>0.
\end{equation}
We use the Taylor expansion of the function $h_k^{(1)}$ around $0$,
$$h_k^{(1)}(\phi_1)=h_k^{(1)}(0)+(h_k^{(1)})'(0)\phi_1+\frac{1}{2}(h_k^{(1)})''(\xi)\phi_1^2, \quad 0<\xi<\phi_1.$$
Set $M_1=\max\vert(h_k^{(1)})''(\xi)\vert<\infty$. Then we have
\begin{equation}\label{lm:taylor1}
h_k^{(1)}(\phi_1)-h_k^{(1)}(0) \geq (h_k^{(1)})'(0)\phi_1-\frac{1}{2}M_1\phi_1^2.
\end{equation}
The derivative of $h_k^{(1)}$ is given by
$$(h_k^{(1)})'(\phi_1)=
\langle\nabla_U f(U_kR(i_k,j_k,\phi_1),V_k,W_k), U_k\dot{R}(i_k,j_k,\phi_1)\rangle.$$
In particular,
\begin{equation}\label{lm:h1-0}
(h_k^{(1)})'(0)=\langle\nabla_U f(U_k,V_k,W_k), U_k\dot{R}(i_k,j_k,0)\rangle.
\end{equation}
It follows from Lemma~\ref{tm:lemma2}(i) and relation~\eqref{lm:h1-0} that
\begin{equation}\label{lm:h1ineq}
\vert(h_k^{(1)})'(0)\vert\geq\eta\|\nabla_U f(U_k,V_k,W_k)\|_2.
\end{equation}
Therefore, from~\eqref{lm:h1ineq} and~\eqref{lm:mu1} we get
\begin{equation}\label{lm:etamu1}
\vert(h_k^{(1)})'(0)\vert\geq\eta\mu_1>0.
\end{equation}
We go back to the inequality~\eqref{lm:taylor1}. For $\phi_1=\frac{1}{M_1}(h_k^{(1)})'(0)$, using the definition of the function $h_k^{(1)}$ and relations~\eqref{lm:fineq} and~\eqref{lm:etamu1}, we obtain
\begin{align*}
& f(U_{k+1},V_{k+1},W_{k+1})-f(U_k,V_k,W_k) \\
& \geq f(U_{k+1},V_k,W_k)-f(U_k,V_k,W_k) =h_k^{(1)}(\phi_1)-h_k^{(1)}(0) \geq(h_k^{(1)})'(0)\phi_1-\frac{1}{2}M_1\phi_1^2 \\ & =\frac{1}{M_1}((h_k^{(1)})'(0))^2-\frac{1}{2M_1}((h_k^{(1)})'(0))^2 \geq \frac{\eta^2\mu_1^2}{2M_1}=\delta>0.
\end{align*}

Now we assume that $\nabla_U f(\overline{U},\overline{V},\overline{W})=0$ and $\nabla_V f(\overline{U},\overline{V},\overline{W})\neq0$.
There is an $\epsilon>0$ such that
\begin{equation}\label{lm:mu2}
\mu_2:=\min\{\|\nabla_V f(U,V,W)\|_2 \ : \ \|V-\overline{V}\|_2<\epsilon\}>0.
\end{equation}
In this case we use the Taylor expansion of the function $h_k^{(2)}$ around $0$. We have
\begin{equation}\label{lm:taylor2}
h_k^{(2)}(\phi_2)-h_k^{(2)}(0) \geq (h_k^{(2)})'(0)\phi_2-\frac{1}{2}M_2\phi_1^2,
\end{equation}
for $M_2=\max\vert(h_k^{(2)})''(\xi)\vert<\infty$, and
$$(h_k^{(2)})'(0)=\langle\nabla_V f(U_k,V_k,W_k), V_k\dot{R}(i_k,j_k,0)\rangle.$$
Lemma~\ref{tm:lemma2}(ii) and relation~\eqref{lm:mu2} imply
\begin{equation}\label{lm:etamu2}
\vert(h_k^{(2)})'(0)\vert\geq\eta\|\nabla_V f(U_k,V_k,W_k)\|_2\geq\eta\mu_2>0.
\end{equation}
The assertion of the lemma follows from~\eqref{lm:taylor2}, \eqref{lm:fineq}, and~\eqref{lm:etamu2} with $\phi_2=\frac{1}{M_2}(h_k^{(2)})'(0)$.

Finally, if $\nabla_U f(\overline{U},\overline{V},\overline{W})=0$ and $\nabla_V f(\overline{U},\overline{V},\overline{W})=0$,
since $\nabla f(U,V,W)(\overline{U},\overline{V},\overline{W})\neq0$,
then it must be that $\nabla_W f(\overline{U},\overline{V},\overline{W})\neq0$.
Then, there is an $\epsilon>0$ such that
\begin{equation}\label{lm:mu3}
\mu_3:=\min\{\|\nabla_W f(U,V,W)\|_2 \ : \ \|W-\overline{W}\|_2<\epsilon\}>0.
\end{equation}
Here we need the Taylor expansion of $h_k^{(3)}$ around $0$,
\begin{equation}\label{lm:taylor3}
h_k^{(3)}(\phi_3)-h_k^{(3)}(0) \geq (h_k^{(3)})'(0)\phi_1-\frac{1}{2}M_3\phi_3^2,
\end{equation}
for $M_3=\max\vert(h_k^{(2)})''(\xi)\vert<\infty$. We repeat the same steps as for the preceding two cases. We have
$$(h_k^{(3)})'(0)=\langle\nabla_W f(U_k,V_k,W_k), W_k\dot{R}(i_k,j_k,0)\rangle,$$
and, using Lemma~\ref{tm:lemma2}(iii) and the relation~\eqref{lm:mu3}, it follows that
\begin{equation}\label{lm:etamu3}
\vert(h_k^{(3)})'(0)\vert\geq\eta\|\nabla_W f(U_k,V_k,W_k)\|_2\geq\eta\mu_3>0.
\end{equation}
We attain inequality~\eqref{lm:fdelta} using~\eqref{lm:taylor3}, \eqref{lm:fineq}, and~\eqref{lm:etamu3} with $\phi_3=\frac{1}{M_3}(h_k^{(3)})'(0)$.
\end{proof}

Using Lemma~\ref{tm:lemma3} we can now prove Theorem~\ref{tm:cvg}.

\begin{proof}[Proof of Theorem~\ref{tm:cvg}]
Suppose that $\overline{U}$, $\overline{V}$, $\overline{W}$ are, respectively, accumulation points of the sequences $\{U_j\}_{j\geq1}$, $\{V_j\}_{j\geq1}$, $\{W_j\}_{j\geq1}$ generated by Algorithm~\ref{agm:jacobi}.
Then there are subsequences $\{U_j\}_{j\in K_U}$, $\{V_j\}_{j\in K_V}$, $\{W_j\}_{j\in K_W}$  such that
$$\{U_j\}_{j\in K_U}\rightarrow\overline{U}, \quad \{V_j\}_{j\in K_V}\rightarrow\overline{V}, \quad \{W_j\}_{j\in K_W}\rightarrow\overline{W},$$
where $K_U,K_V,K_W\subseteq\mathbb{N}$.

Assume that $(\overline{U},\overline{V},\overline{W})$ is not a stationary point of the function $f$, that is
\begin{equation}\label{tm:contradiction}
\nabla f(\overline{U},\overline{V},\overline{W})\neq0.
\end{equation}
Then, for any $\epsilon>0,$ there are $k_0^{(U)}\in K_U$, $k_0^{(V)}\in K_V$, $k_0^{(W)}\in K_W$ such that
$$\|U_k-\hat{U}\|_2<\epsilon, \quad \|V_k-\hat{V}\|_2<\epsilon, \quad \|W_k-\hat{W}\|_2<\epsilon,$$
for every $k>k_0$, $k_0=\max\{k_0^{(U)},k_0^{(V)},k_0^{(W)}\}$.
Thus, Lemma~\ref{tm:lemma3} implies $f(U_{k+1},V_{k+1},W_{k+1})-f(U_k,V_k,W_k)\geq\delta>0$.
It follows that
$$f(U_k,V_k,W_k)\rightarrow\infty,$$
when $k\rightarrow\infty$.
Since $f$ is continuous, if $(U_k,V_k,W_k)$ converges, then $f(U_k,V_k,W_k)$ should converge, too. This gives a contradiction. Therefore, assumption~\eqref{tm:contradiction} cannot hold and $(\overline{U},\overline{V},\overline{W})$ is a stationary point of $f$.
\end{proof}

Note that all results from this section can be generalized to order-$d$ tensors, $d>3$.

\section{Numerical examples}\label{sec:numerical}

We illustrate the convergence of Algorithm~\ref{agm:jacobi} through several numerical examples.
We observe the relative off-norm of a tensor $\calA$, which is given as
\begin{equation}\label{num:reloff}
\frac{\off(\calA)}{\|\calA\|}.
\end{equation}
For a diagonal tensor, value~\eqref{num:reloff} is equal to zero, while for a random tensor it is, typically, close to one. Note that the off-norm is not a norm because it can be equal to zero for a nonzero input.

\begin{figure}[h!]
    \centering
    \includegraphics[width=0.45\textwidth]{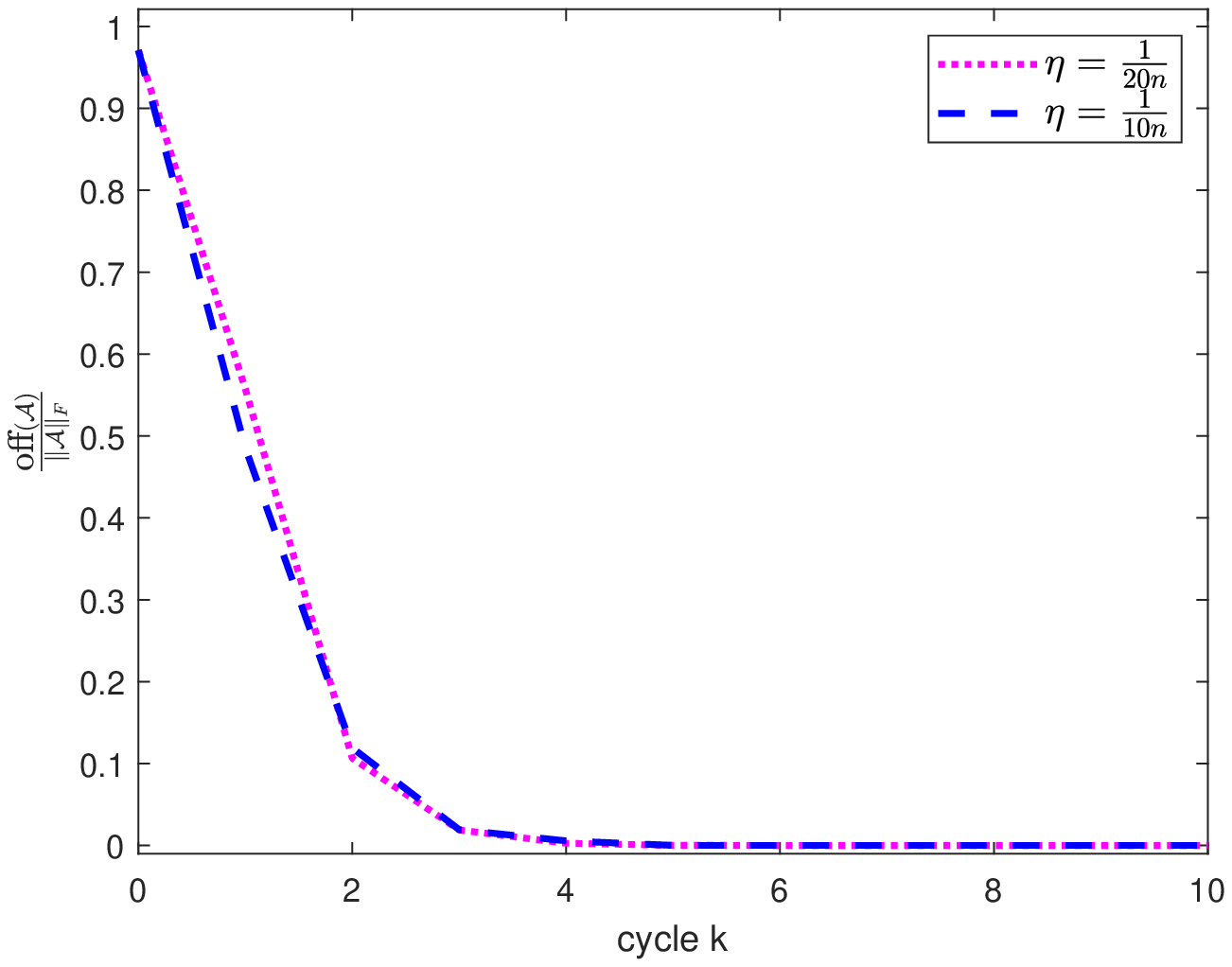}
    \includegraphics[width=0.45\textwidth]{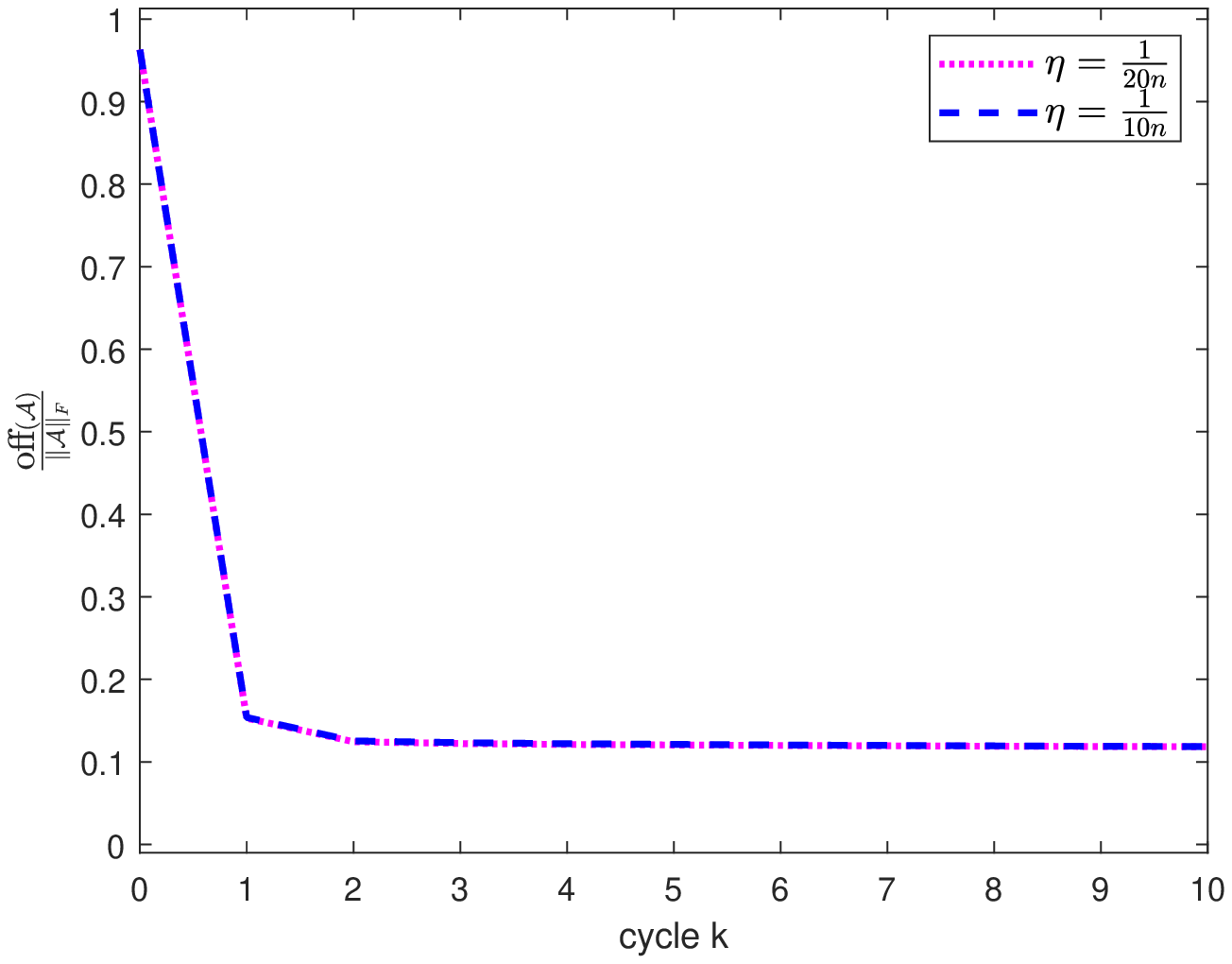}
    \caption{Change in the relative off-norm for two $30\times30\times30$ tensors with different values of $\eta$. Left: Diagonalizable tensor. Right: Non-diagonalizable tensor.}\label{fig:test0}
\end{figure}

Figure~\ref{fig:test0} shows the change of the relative off-norm. We distinguish two different situations, one where a tensor can be completely diagonalized using orthogonal transformations, and a more general one where orthogonal diagonalization is not possible. For the first set of tensors we get $\frac{\off(\calA)}{\|\calA\|}=0$. Otherwise, we get the convergence to some value between $0$ and $1$. A random diagonalizable tensor is constructed by taking a diagonal tensor with random entries on the diagonal and multiplying it in each mode by orthogonal matrices obtained from QR factorizations of three random matrices. The algorithm uses row-wise cyclic pivot ordering~\eqref{Or} with different values of the parameter $\eta$.

\begin{figure}[h!]
    \centering
    \includegraphics[width=0.45\textwidth]{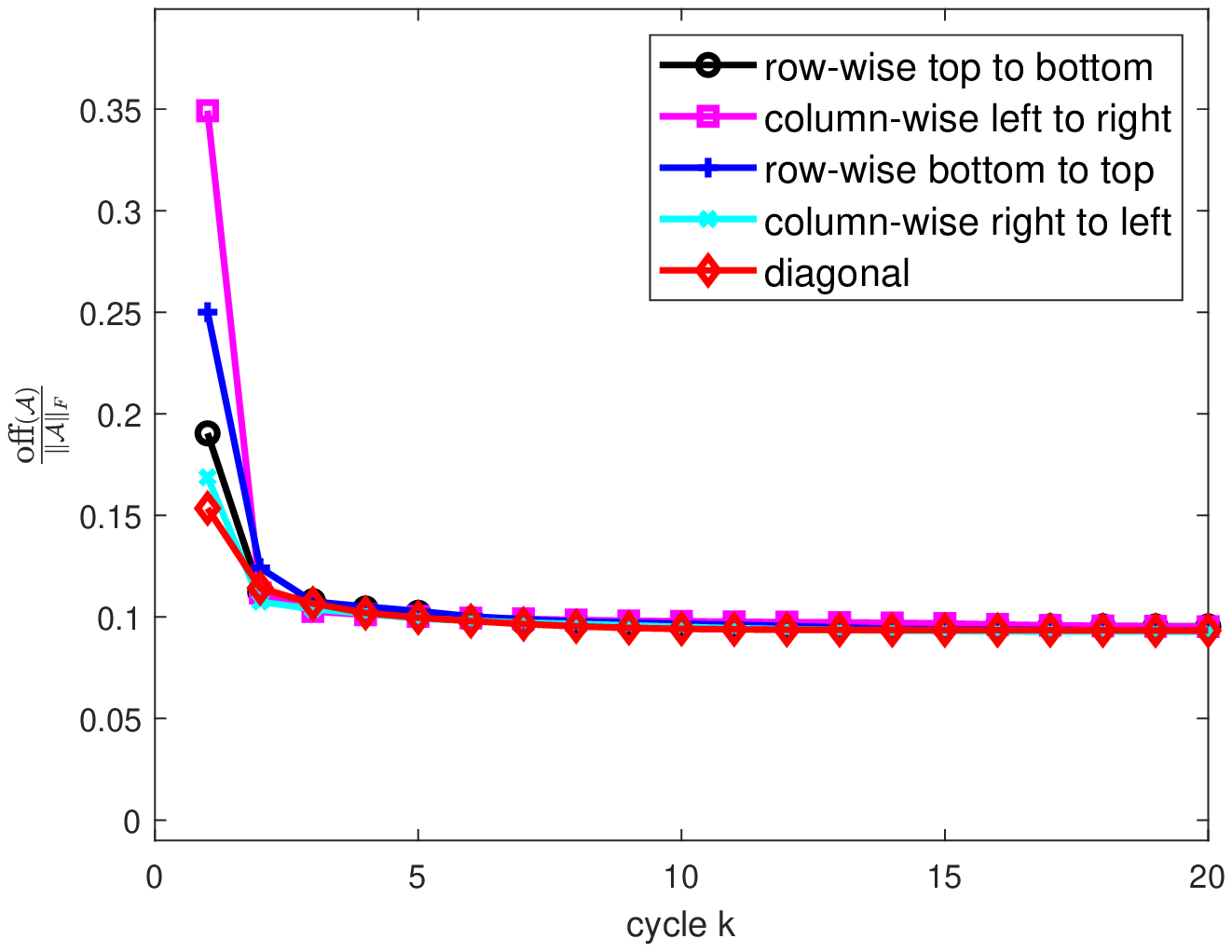}
    \includegraphics[width=0.45\textwidth]{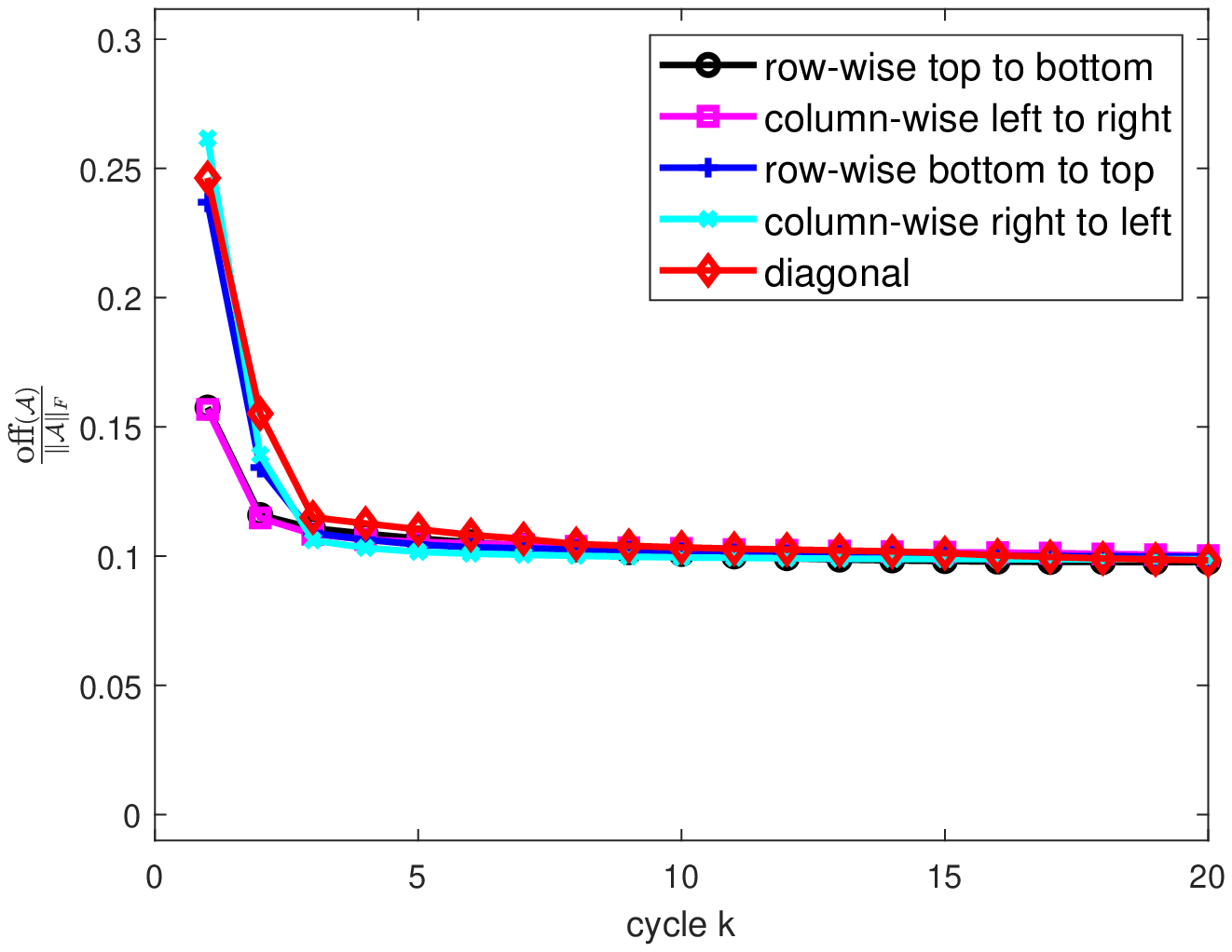}
    \includegraphics[width=0.45\textwidth]{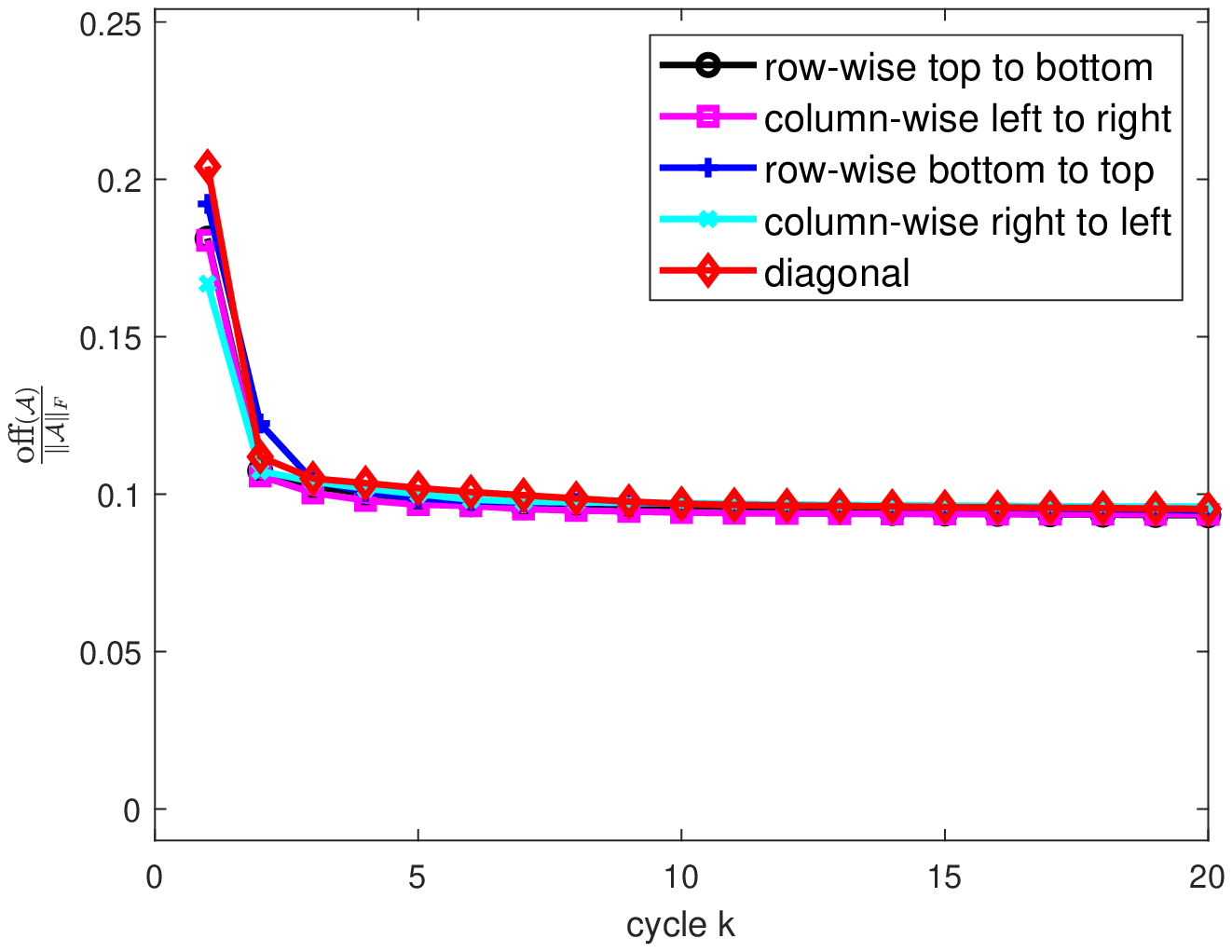}
    \includegraphics[width=0.45\textwidth]{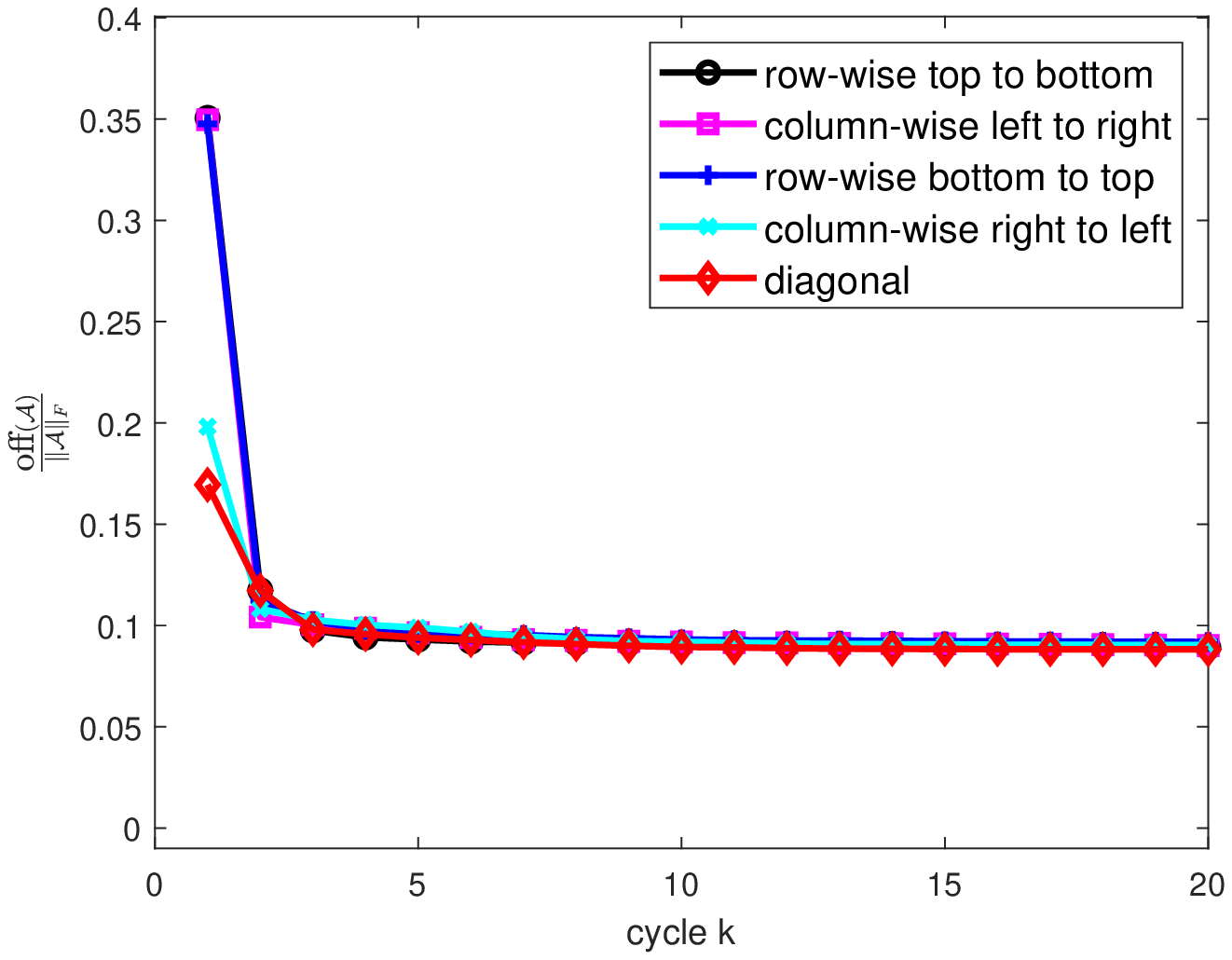}
    \includegraphics[width=0.45\textwidth]{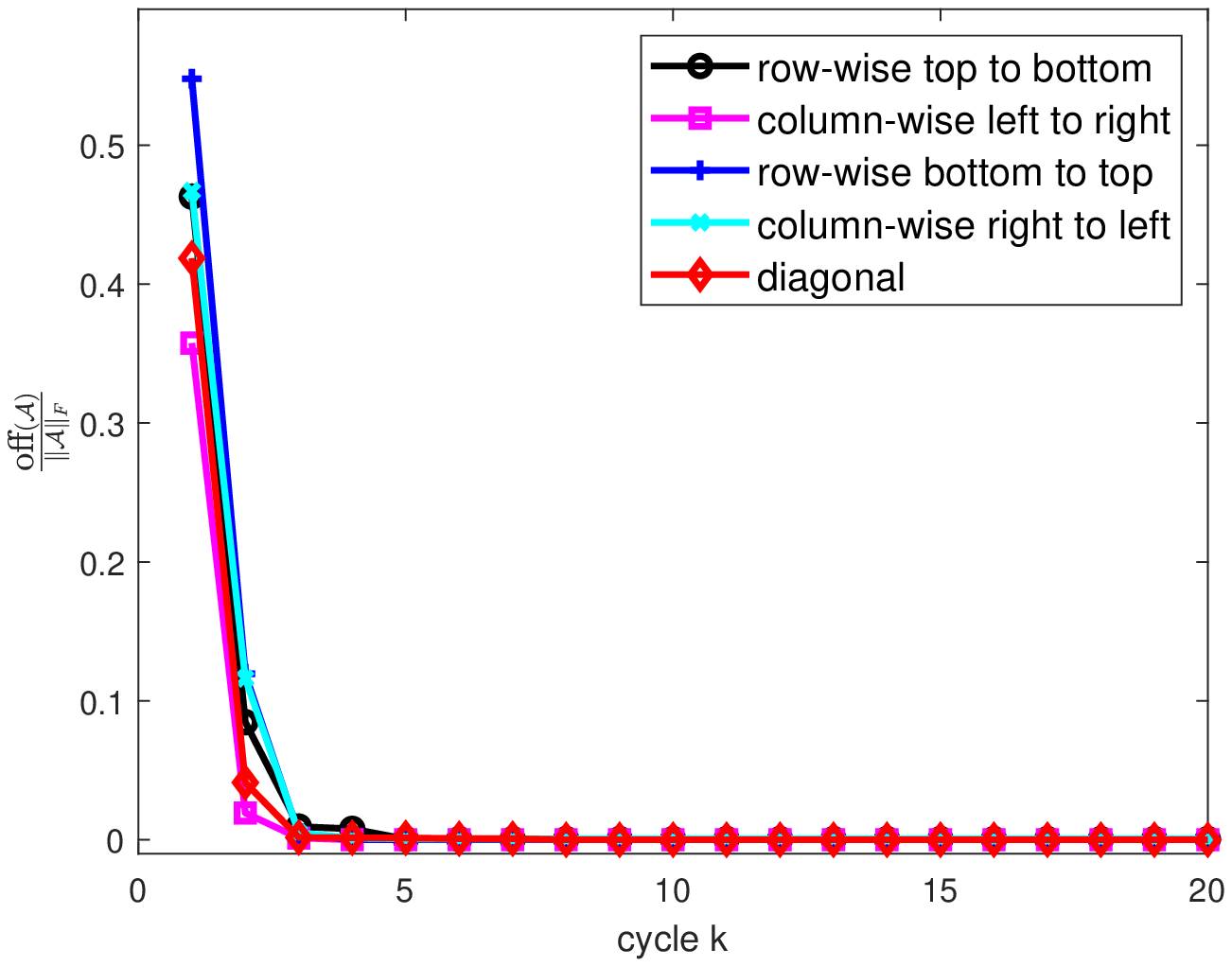}
    \includegraphics[width=0.45\textwidth]{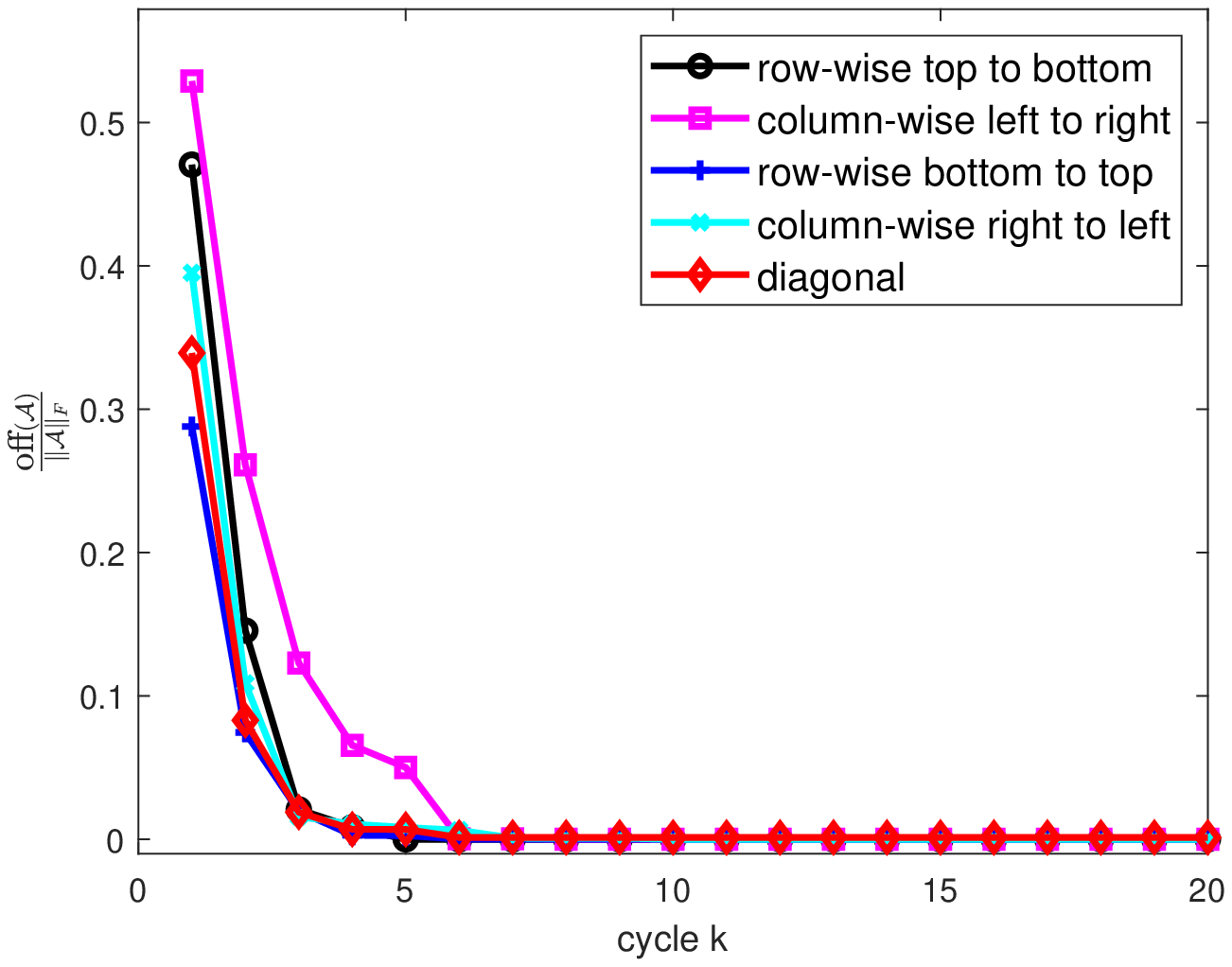}
    \caption{Change in the relative off-norm for six $10\times10\times10$ tensors with different pivot strategies.}\label{fig:O}
\end{figure}

In Figure~\ref{fig:O} we compare five different pivot orderings. In addition to the row-wise top to bottom~\eqref{Or} and the column-wise left to right~\eqref{Oc} ordering we have the row-wise bottom to top
$$\mathcal{O}_r'=(n-1,n),(n-2,n-1),(n-2,n),(n-3,n-2),\ldots,(2,n),(1,2),\ldots,(1,n),$$
the column-wise right to left
$$\mathcal{O}_c'=(1,n),\ldots,(n-1,n),(1,n-1),\ldots,(n-2,n-1),\ldots,(1,2),$$
and the diagonal ordering of pivot pairs
$$\mathcal{O}_d=(1,2),(2,3),(3,4),\ldots,(n-1,n),(1,3),(2,4),\ldots,(n-2,n),(1,4),\ldots,(1,n).$$
We run the algorithm on six random tensors, four of which cannot be diagonalized by orthogonal transformations, all with $\eta=\frac{1}{20n}$. As expected, different pivot strategies are faster/slower on different tensors. However, no matter what pivot strategy we choose, the algorithm converges to the same point.

\begin{figure}[h!]
    \centering
    \includegraphics[width=0.45\textwidth]{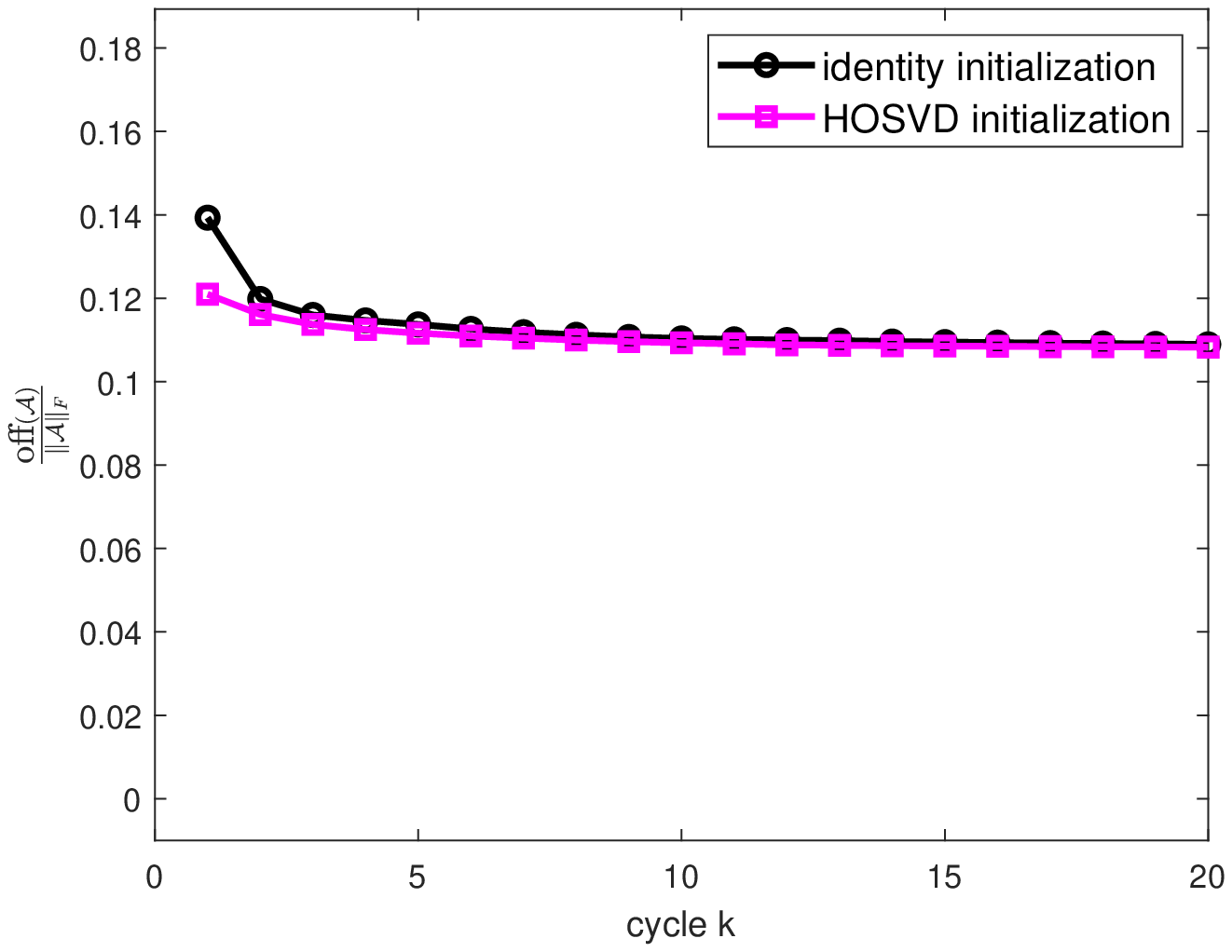}
    \includegraphics[width=0.45\textwidth]{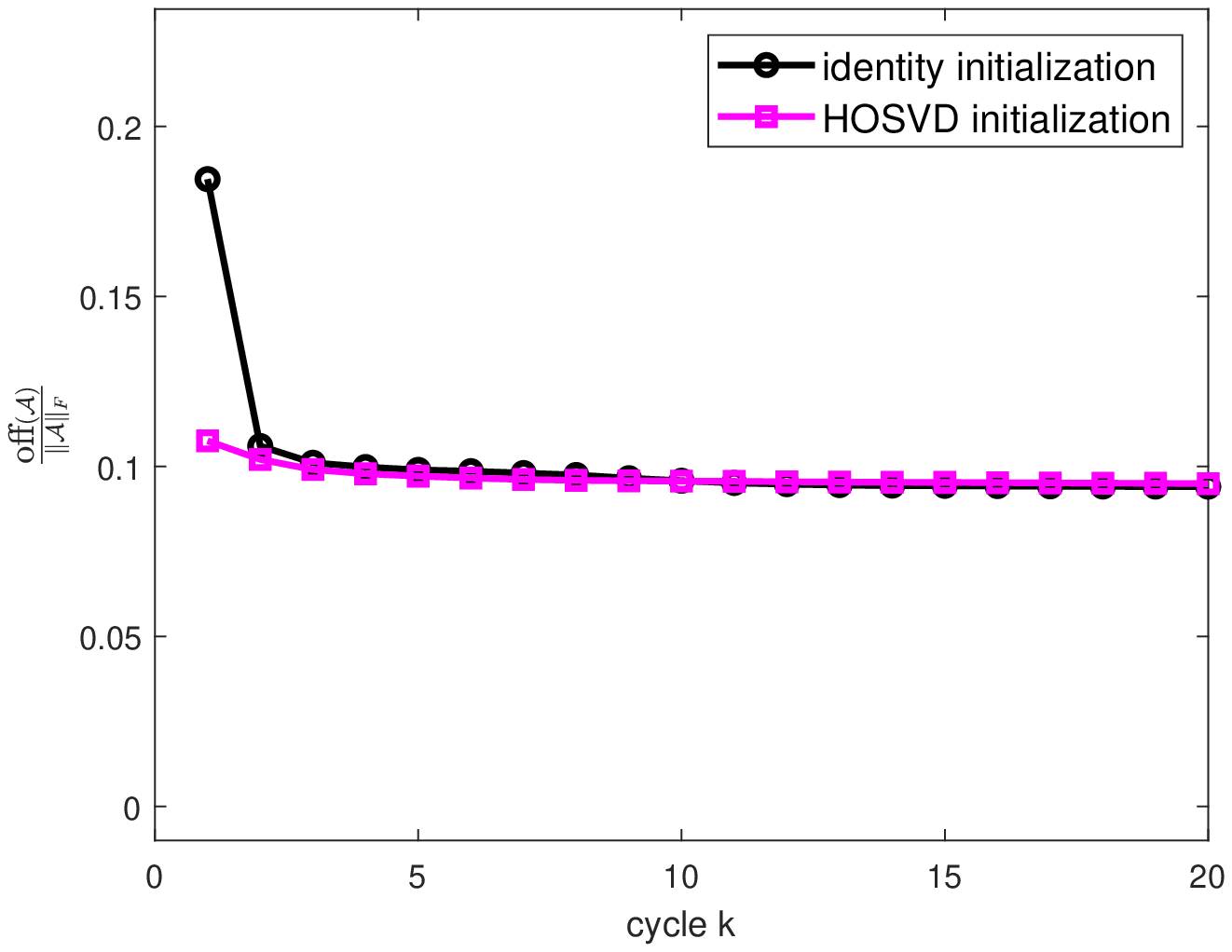}
    \caption{Convergence of the algorithm with different initializations. Left: $20\times20\times20$ tensor. Right: $10\times10\times10$ tensor.}\label{fig:test1}
\end{figure}

In Figure~\ref{fig:test1} we compare two different initializations for our Jacobi-type algorithm. The first one is identity initialization as it was done in Algorithm~\ref{agm:jacobi}. There we have
\begin{equation}\label{num:id}
\calA^{(0)}=\calA, \quad U_0=V_0=W_0=I.
\end{equation}
The other initialization that can be used is coming from the HOSVD of $\calA$, see~\cite{DeLHOSVD},
$\calA=\widetilde{\calS}\times_1\widetilde{U}\times_2\widetilde{V}\times_3\widetilde{W},$
where $\widetilde{U}$, $\widetilde{V}$, and $\widetilde{W}$ are matrices of left singular vectors of matricizations $A_{(1)}$, $A_{(2)}$, and $A_{(3)}$, respectively, and
$\widetilde{\calS}=\calA\times_1\widetilde{U}^T\times_2\widetilde{V}^T\times_3\widetilde{W}^T.$
Then, instead of the initialization~\eqref{num:id}, we set
$\calA^{(0)}=\widetilde{\calS}$, $U_0=\widetilde{U}$, $V_0=\widetilde{V}$, $W_0=\widetilde{W}$.
We run the algorithm with $\eta=\frac{1}{20n}$ on two random tensors. We can see that the HOSVD initialization is superior in the beginning cycles. This is the case because, compared to the starting tensor $\calA$, the core tensor $\widetilde{\calS}$ from the HOSVD of $\calA$ is significantly closer to a diagonal tensor. Nevertheless, after those first cycles, both initializations are equally good.

\section{Symmetric and antisymmetric tensors}\label{sec:structured}

We say that a tensor is symmetric if its elements remain constant under any permutation of indices. For a symmetric tensor $\calX\in\mathbb{R}^{n\times n\times n}$ we have
$$x_{ijk}=x_{ikj}=x_{jik}=x_{jki}=x_{kij}=x_{kji}.$$
Symmetric tensors were studied in details in~\cite{LUC18} and~\cite{LUC19}, where the authors also work with a Jacobi-type algorithm.

Our algorithm is not structure-preserving. In order to have a symmetry-preserving Jacobi-type algorithm, rotation matrices should be the same in all modes. Since the rotations $R_{U,k}$, $R_{V,k}$, and $R_{W,k}$ in the $k$th step are chosen depending on different tensors, they are not all the same. Nevertheless, we have noticed in practice that for smaller values of $\eta$ from~\eqref{pivotcond}, after the convergence criterion is satisfied, the Algorithm~\ref{agm:jacobi}, in most of the cases, returns mutually equal matrices $U,V,W$ and a symmetric tensor $\calS$. However, this is not the case for larger $\eta$. We will illustrate this behaviour on an example.

\begin{figure}[h!]
    \centering
    \includegraphics[width=0.45\textwidth]{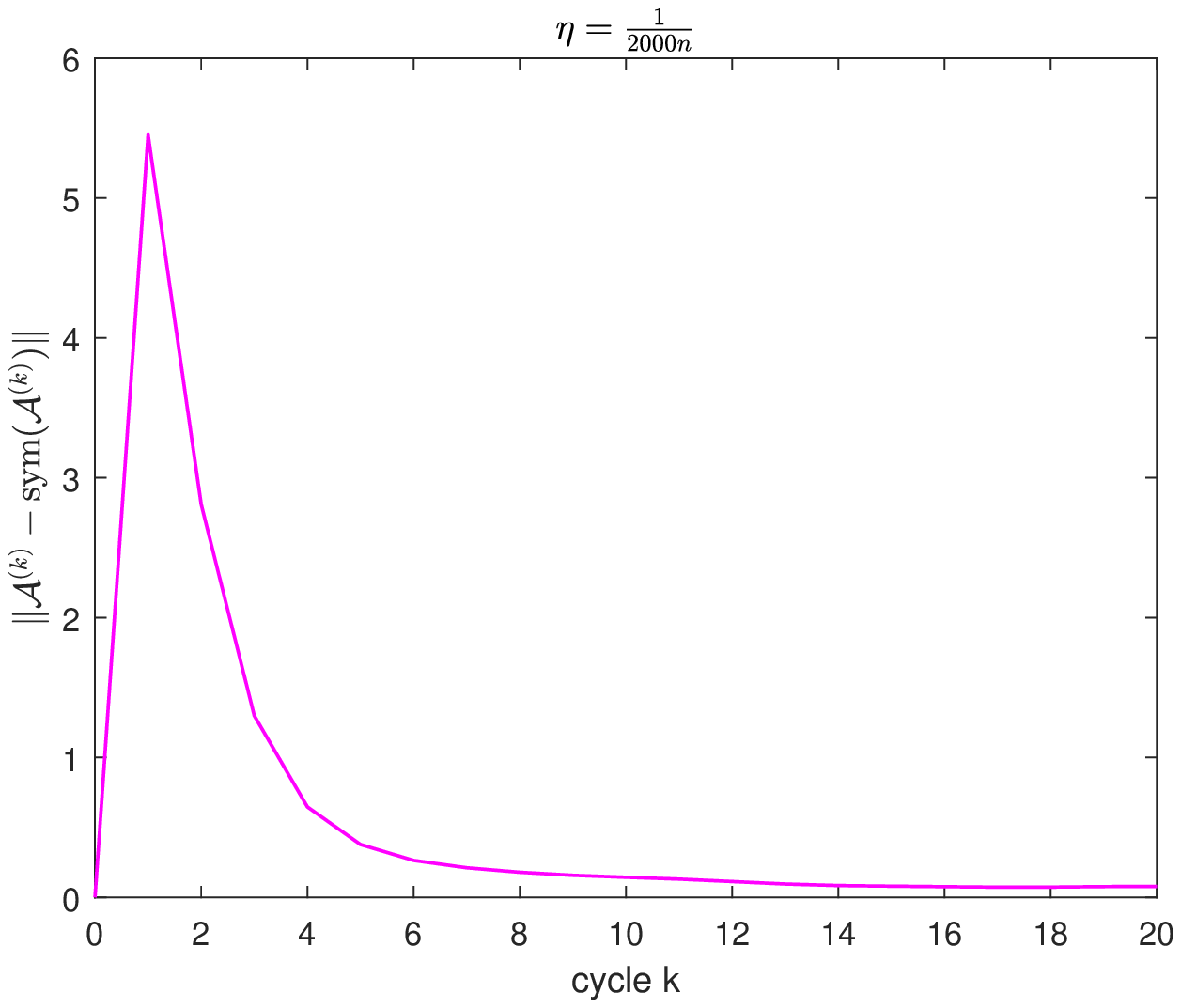}
    \includegraphics[width=0.45\textwidth]{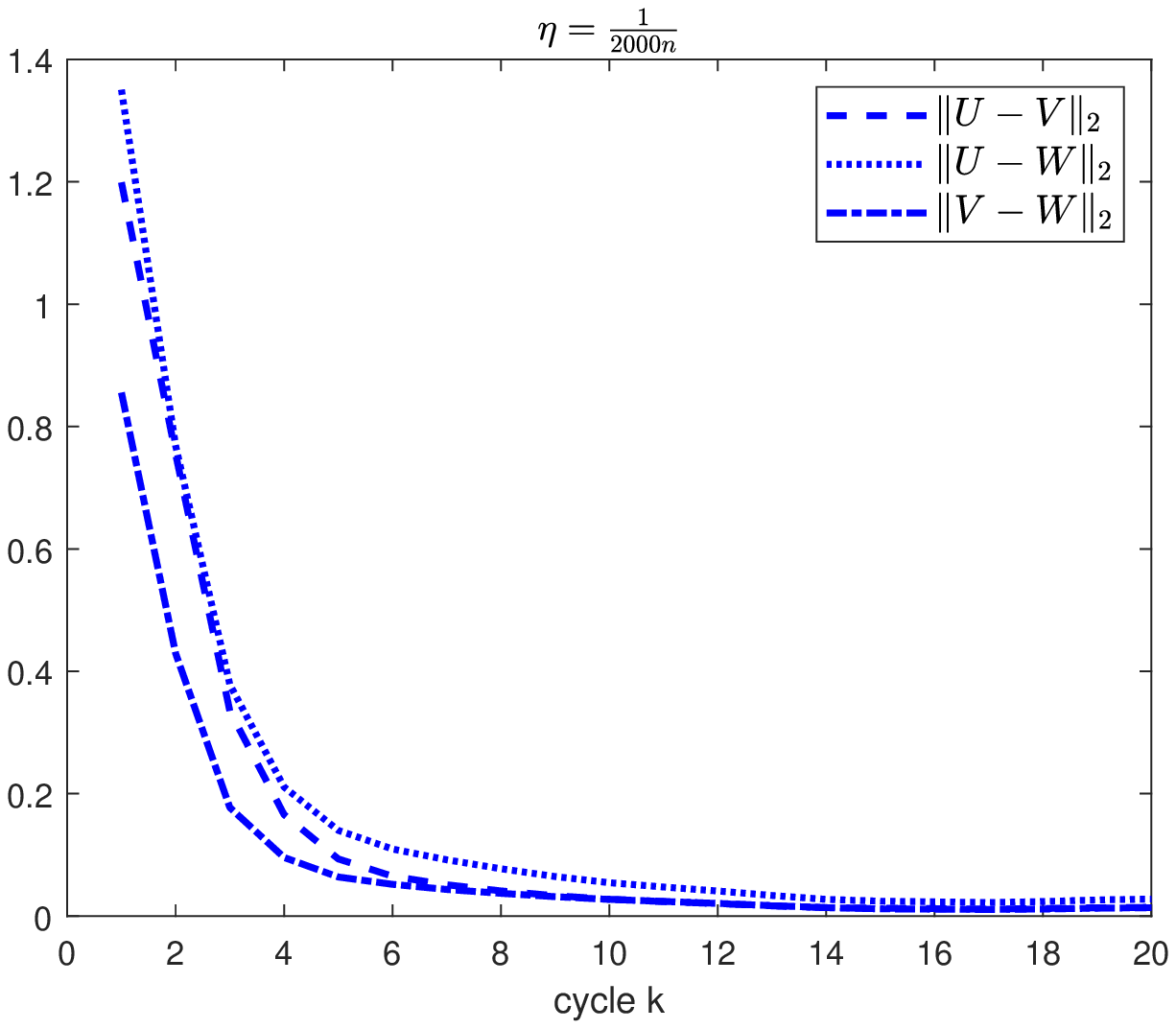}
    \caption{Departure from the symmetry for a random symmetric $20\times20\times20$ tensor.}\label{fig:distancesym}
\end{figure}

Departure from symmetry is measured as the distance in the Frobenius norm between the tensor $\calA$ and its symmetrization $\text{sym}(\calA)$,
\begin{equation}\label{symdistance}
\|\calA-\text{sym}(\calA)\|,
\end{equation}
where, for a $3$rd order tensor $\calX$, we have
$$\text{sym}(\calX)=\frac{1}{6}(x_{ijk}+x_{ikj}+x_{jik}+x_{jki}+x_{kij}+x_{kji}).$$
It is easy to check that, if $\calX$ is symmetric, then $\text{sym}(\calX)=\calX$,
and the expression~\eqref{symdistance} is equal to zero.
We applied Algorithm~\ref{agm:jacobi} with $\eta=\frac{1}{2000n}$ on a randomly generated symmetric $20\times20\times20$ tensor $\calA$. In the left picture in Figure~\ref{fig:distancesym} we can see that after we start with a symmetric tensor, symmetry is lost already in the first cycle, as expected, but the tensor becomes more and more symmetric through iterations and the sequence $(\calA^{(k)})_k$ converges to a symmetric tensor. In the right picture we see that the distance between each pair of matrices $U_k,V_k,W_k$ converges to zero, that is, $U_k,V_k,W_k$ converge to the same matrix. This does not happen for $\eta=\frac{1}{20n}$.

One should keep in mind that for a symmetric starting tensor the solution of the maximization problem~\eqref{Smaximization} is not necessary a symmetric tensor. One such example is tensor $\mathcal{T}$ from~\cite[Example 5.5]{LUC19} that can be given by its matricization
\begin{equation}\label{tensorCom}
T_{(1)}={\footnotesize\left[
                \begin{array}{ccccccccc}
                  0 & 0 & 0 & 0 & 0 & 1 & 0 & 1 & 0 \\
                  0 & 0 & 1 & 0 & 0 & 0 & 1 & 0 & 0 \\
                  0 & 1 & 0 & 1 & 0 & 0 & 0 & 0 & 0 \\
                \end{array}
              \right]}.
\end{equation}
First, we notice that neither identity nor HOSVD initialization work on this tensor. In both cases the diagonal elements of $\mathcal{T}^{(0)}$ are zero and all rotation angles are zero, so the tensor is unchanged. Thus, here we do preconditioning with a random orthogonal matrix $Q$ by setting
$\mathcal{A}^{(0)}=\mathcal{A}\times_1Q\times_2Q\times_3Q.$
For the vast majority of the choices of $Q$, Algorithm~\ref{agm:jacobi} with $\eta=\frac{1}{2000n}$ converges to the symmetric tensor $\mathcal{S}$,
$$S_{(1)}={\footnotesize\left[
                \begin{array}{ccccccccc}
                  0.8889 & -0.4444 & -0.4444 & -0.4444 & -0.4444 & -0.1111 & -0.4444 & -0.1111 & -0.4444 \\
                  -0.4444 & -0.4444 & -0.1111 & -0.4444 & 0.8889 & -0.4444 & -0.1111 & -0.4444 & -0.4444 \\
                  -0.4444 & -0.1111 & -0.4444 & -0.1111 & -0.4444 & -0.4444 & -0.4444 & -0.4444 & 0.8889 \\
                \end{array}
              \right]},
$$
with transformation matrix $U=V=W$ depending on $Q$.
This is a stationary point of the objective function~\eqref{maxfunction}, but not a point of its global maximum.
In the other rare cases the algorithm converged to one of the better, but nonsymmetric, solutions of the form $\bar{\mathcal{S}}$
$$\bar{S}_{(1)}={\footnotesize\left[
                \begin{array}{ccccccccc}
                  \pm1 & 0 & 0 & 0 & 0 & 0 & 0 & \pm1 & 0 \\
                  0 & 0 & \pm1 & 0 & \pm1 & 0 & 0 & 0 & 0 \\
                  0 & 0 & 0 & \pm1 & 0 & 0 & 0 & 0 & \pm1 \\
                \end{array}
              \right]}.
$$

On the other hand, a tensor is antisymmetric if its elements change sign when permuting pairs of indices. For an antisymmetric tensor $\calX\in\mathbb{R}^{n\times n\times n}$ we have
$$x_{ijk}=x_{jki}=x_{kij}=-x_{ikj}=-x_{jik}=-x_{kji}.$$
In every antisymmetric tensor, elements on the positions where two or more indices are the same are equal to zero. Hence, all elements on the diagonal of an antisymmetric tensor are zero. This is the reason why, contrary to the symmetric case where one may be interested in a structure-preserving algorithm, we are not interested in preserving the antisymmetry.
Here, by each iteration a tensor moves further from the structure.
Still, antisymmetric tensors need some special attention.
If we apply the algorithm directly in the form given in Algorithm~\ref{agm:jacobi}, with the identity initialization, the algorithm will fail when computing the rotation angle. This happens because for an antisymmetric tensor, when computing the tangent of the double rotation angle~\eqref{tan2phi1}, \eqref{tan2phi2}, and~\eqref{tan2phi3}, we get both the numerator and the denominator equal to zero. We can overcome this problem with a preconditioning step --- instead of the identity initialization~\eqref{num:id} we use the HOSVD initialization as described in Section~\ref{sec:numerical}.

\end{document}